\documentclass[12pt]{amsart}
\usepackage{tikz}
\usepackage{tikz-cd}
\usetikzlibrary{matrix, calc, arrows,backgrounds}

\tikzset{node distance=2cm,auto}
\usepackage{amsfonts,amsmath,amssymb,amsthm,mathrsfs,verbatim}
\usepackage[all,cmtip]{xy}
\usepackage{latexsym}
\usepackage{graphics}
\newtheorem{theorem}{Theorem}[section]
\newtheorem{lemma}[theorem]{Lemma}
\newtheorem{corollary}[theorem]{Corollary}
\newtheorem{proposition}[theorem]{Proposition}
\newtheorem{conjecture}{Conjecture}

\newtheorem{mainthm}{Theorem}

\theoremstyle{definition}
\newtheorem{question}[theorem]{Question}
\newtheorem{definition}[theorem]{Definition}
\newtheorem{example}[theorem]{Example}

\newtheorem{remark}[theorem]{Remark}

\usepackage{amsfonts}
\usepackage{color}

\newcommand{\bpf}{\noindent{\bf Proof}\hspace{7pt}}
\newcommand{\epf}{\qed}
\newcommand{\ben}{\begin{enumerate}}
\newcommand{\een}{\end{enumerate}}
\newcommand{\ble}{\begin{lemma}}
\newcommand{\ele}{\end{lemma}}
\newcommand{\bth}{\begin{theorem}}
\renewcommand{\eth}{\end{theorem}}
\newcommand{\bmth}{\begin{mainthm}}
\newcommand{\emth}{\end{mainthm}}

\newcommand{\bpr}{\begin{proposition}}
\newcommand{\epr}{\end{proposition}}
\newcommand{\bco}{\begin{corollary}}
\newcommand{\eco}{\end{corollary}}
\newcommand{\bcon}{\begin{conjecture}}
\newcommand{\econ}{\end{conjecture}}
\newcommand{\bqu}{\begin{question}}
\newcommand{\equ}{\end{question}}
\newcommand{\bde}{\begin{definition}}
\newcommand{\ede}{\end{definition}}
\newcommand{\bre}{\begin{remark}}
\newcommand{\ere}{\end{remark}}
\newcommand{\bex}{\begin{example}}
\newcommand{\eex}{\end{example}}
\newcommand{\barr}{\begin{array}}
\newcommand{\earr}{\end{array}}
\newcommand{\btab}{\begin{tabular}}
\newcommand{\etab}{\end{tabular}}
\newcommand{\beq}{\begin{equation}}
\newcommand{\eeq}{\end{equation}}
\newcommand{\bea}{\begin{eqnarray*}}
\newcommand{\eea}{\end{eqnarray*}}
\newcommand{\bce}{\begin{center}}
\newcommand{\ece}{\end{center}}
\newcommand{\bpi}{\begin{picture}}
\newcommand{\epi}{\end{picture}}
\newcommand{\bfi}{\begin{figure} \begin{center}}
\newcommand{\efi}{\end{center} \end{figure}}

\newcommand{\bsl}{\begin{slide}{}}
\newcommand{\esl}{\end{slide}}

\newcommand{\pf}{\noindent{\bf Proof}\hspace{7pt}}
\newcommand{\ul}{\underline}

\newcommand{\hso}[1]{\hspace{-1pt}}

\newcommand{\into}{\hookrightarrow}
\newcommand{\onto}{\twoheadrightarrow}

\newcommand{\sbe}{\subseteq}

\def\<{\langle}
\def\>{\rangle}

\newcommand{\De}{\Delta}

\newcommand{\La}{\Lambda}

\DeclareMathOperator{\edg}{E}

\newcommand{\liediag}{\Gamma}
\newcommand{\covdiag}{\Lambda}
\newcommand{\uliediag}{\ul{\liediag}}
\newcommand{\trin}{\tau}

\newcommand{\cF}{{\mathcal F}}

\newcommand{\cR}{{\mathcal R}}

\newcommand{\sF}{{\mathsf{F}}}

\newcommand{\tA}{\widetilde{A}}

\newcommand{\hI}{{I(\La)}}

\newcommand{\amgrp}[1]{{\mathbf{#1}}}
\newcommand{\comp}[1]{{\mathrm{#1}}}

\newcommand{\amgrpC}{{\mathbf{C}}}

\newcommand{\amgrpB}{{\mathbf{B}}}
\newcommand{\amgrpF}{{\mathbf{F}}}
\newcommand{\amgrpG}{{\mathbf{G}}}
\newcommand{\amgrpH}{{\mathbf{H}}}
\newcommand{\amgrpL}{{\mathbf{L}}}

\newcommand{\amgrpX}{{\mathbf{X}}}

\newcommand{\amf}{{\mathbf f}}
\newcommand{\famf}{\ul{\mathbf f}}
\newcommand{\amg}{{\mathbf g}}
\newcommand{\famg}{\ul{\mathbf g}}

\newcommand{\aml}{{\mathbf l}}
\newcommand{\uaml}{\tilde{\mathbf l}}

\newcommand{\bamgrpF}{\bar{\amgrpF}}
\newcommand{\bamgrpG}{\bar{\amgrpG}}

\newcommand{\uammapp}{\tilde{\pi}}

\newcommand{\compD}{{{D}}}
\newcommand{\compF}{{{F}}}
\newcommand{\compG}{{{G}}}

\newcommand{\compL}{{{L}}}
\newcommand{\compN}{{{N}}}

\newcommand{\compB}{{{B}}}
\newcommand{\compU}{{{U}}}

\newcommand{\compl}{\lambda}
\newcommand{\compf}{\phi}
\newcommand{\compg}{\gamma}

\newcommand{\ucompG}{{\tilde{G}}}

\newcommand{\ucompL}{{\tilde{L}}}

\newcommand{\ucompg}{{\tilde{\gamma}}}
\newcommand{\ucompl}{{\tilde{\lambda}}}

\newcommand{\umap}{\pi}

\newcommand{\twA}{{}^2\! {A}}

\newcommand{\Aut}{\mathop{\rm Aut}\nolimits}

\renewcommand{\bar}{\overline}

\newcommand{\id}{\mathop{\rm id}\nolimits}

\newcommand{\Inv}{\mathop{\rm Inv}\nolimits}

\def\flexbox#1{\mathchoice{\mbox{#1}}{\mbox{#1}}{\mbox{\scriptsize #1}}%
{\mbox{\tiny #1}}}
\def\SL{\mathop{\flexbox{\rm SL}}\nolimits}
\def\PSL{\mathop{\flexbox{\rm PSL}}\nolimits}
\newcommand{\POmega}{\mathop{\flexbox{\rm P}\Omega}\nolimits}

\DeclareMathOperator{\SU}{SU}

\newcommand{\CC}{{\mathbb C}}

\newcommand{\FF}{{\mathbb F}}

\newcommand{\NN}{{\mathbb N}}

\newcommand{\ZZ}{{\mathbb Z}}

\newcommand{\hE}{{E(\La)}}

\newcommand{\hpsi}{\lambda}
\newcommand{\bL}{{\mathbf{L}}}

\newcommand{\Pig}{\pi(\Gamma,{0})}

\newcommand\Gm{\Gamma}

\newcommand{\opp}{\mathbin{\rm opp}}

\newcommand{\fk}{{\rm k}}

\newcommand{\dfn}{\em}

\newcommand{\after}{\mathbin{ \circ }}

\newcommand{\Sp}{\mathop{\rm Sp}\nolimits}

\newcommand{\Stab}{\mathop{\rm Stab}}

\DeclareMathOperator{\proj}{proj}

\newcommand{\Chi}{{\mathcal X}}

\renewcommand{\hat}{\widehat}

\newcommand{\vep}{\varepsilon}
\newcommand{\vph}{\varphi}

\renewcommand{\qed}{\hfill $\square$}

\newcounter{romanlistctr}
{\end{list}}%

 \makeatletter
 \@addtoreset{equation}{section}
 \makeatother

 \makeatletter
 \def\section{\@startsection {section}{1}{\z@}{-1.5ex plus -.5ex
 minus -.2ex}{1ex plus .2ex}{\large\bf}}

 \makeatletter
 \def\subsection{\@startsection {subsection}{1}{\z@}{-1.5ex plus -.5ex
 minus -.2ex}{1ex plus .2ex}{\bf}}

\newcommand{\G}{{\mathbf{G}}}

\def\<{\langle}
\def\>{\rangle}

\newcommand{\amA}{{\mathscr{A}}}

\newcommand{\amF}{{\mathscr{F}}}
\newcommand{\famF}{\ul{\mathscr{F}}}

\newcommand{\amG}{{\mathscr{G}}}

\newcommand{\amL}{{\mathscr{L}}}
\newcommand{\famG}{\ul{\mathscr{G}}}

\newcommand{\bhKM}{\bar{L}}

\newcommand{\classmap}{\omega}
\usepackage{hyperref}
\newcommand{\ead}{\email}
\newenvironment{keyword}{Keywords: \keywords}{}

\newcommand{\MSC}[1]{MSC #1:\subjclass}
\newcommand{\sep}{\hspace{1ex}}

\begin{document}
\title[Curtis-Tits groups and Phan groups]{\Large Realizations and properties of $3$-spherical Curtis-Tits Groups and Phan groups}

\author[R.~Blok]{Rieuwert J. Blok}
\ead{blokr@member.ams.org}
\address{Department of Mathematics and Statistics\\
Bowling Green State University\\
Bowling Green, OH 43403\\
U.S.A.}
\author[C.~Hoffman]{Corneliu G. Hoffman}
\ead{C.G.Hoffman@bham.ac.uk}
\address{University of Birmingham\\
Edgbaston, B15 2TT\\
U.K.}

\begin{abstract}
In this note we establish the existence of all Curtis-Tits groups and Phan groups with $3$-spherical diagram as classified in~\cite{BloHofShp2017} and investigate some of their geometric and group theoretic properties.
Whereas it is known that orientable Curtis-Tits groups with spherical or non-spherical and non-affine diagram are almost simple, we show that non-orientable Curtis-Tits groups are acylindrically hyperbolic and therefore have infinitely many infinite-index normal subgroups.
However, we also provide concrete examples of non-orientable Curtis-Tits groups whose quotients are finite simple groups of Lie type.
 \end{abstract}
\maketitle

\begin{keyword}
Curtis-Tits groups, Phan groups, groups of Kac-Moody type, lattices, abstract simplicity.
\MSC[2010] 20G35 \sep 
51E24%
\end{keyword}

\section{Introduction}
In~\cite{BloHofShp2017} Curtis-Tits amalgams and Phan amalgams over a finite field $\FF_q$ with $3$-spherical diagram $\liediag$ and weak system of fundamental root groups (Curtis-Tits case) or property (D) (Phan case) were completely classified, generalizing the result from~\cite{BloHof2014b} which covered the case of Curtis-Tits amalgams over a field of order $\ge 4$ with simply-laced 
 $3$-spherical diagram satisfying property (D).
In the present paper we shall generalize the results from~\cite{BloHof2016} showing that not only all Curtis-Tits amalgams but also all Phan amalgams classified in~\cite{BloHofShp2017} have non-trivial completions.
More precisely, using the notation of Definition~\ref{dfn:representative amalgams} we have the following.

\begin{mainthm}\label{mainthm:CT realization}
Every Curtis-Tits amalgam of the form $\amG(\delta)$ has a non-trivial completion.
In particular, an arbitrary Curtis-Tits amalgam over $\FF_q$ with $3$-spherical diagram having no $C_2(2)$-subdiagrams has a non-trivial completion if and only if it possesses a weak system of fundamental root groups.
\end{mainthm}
\bpf
The first claim is the content of Theorems~\ref{thm:OCT realization}~and~\ref{thm:NOCT realization}.
We now recall that it was shown in~\cite{BloHofShp2017} that an arbitrary Curtis-Tits amalgam over $\FF_q$ with $3$-spherical diagram possessing some non-trivial completion does have a weak system of fundamental root groups and is isomorphic to $\amG(\delta)$ for some $\delta$, so the second claim follows from the first.
\epf

\bre\label{rem:CT non-universal} 
Note that Curtis-Tits amalgams as defined here are all universal in the sense that the groups appearing in them are universal groups of Lie type. In~\cite{BloHof2016} we also derive existence criteria for Curtis-Tits type amalgams with $3$-spherical simply-laced diagrams that are not universal.
We expect that an analogous treatment of general $3$-spherical Curtis-Tits amalgams will yield a similar result, but we shall not work out the details here.
\ere

\begin{mainthm}\label{mainthm:Phan realization}
Every Phan amalgam of the form $\amG(\delta)$ has a non-trivial completion.
In particular, an arbitrary Phan amalgam over $\FF_q$ with $3$-spherical diagram has a non-trivial completion if and only if it satisfies property (D).
\end{mainthm}
\bpf
The first claim is the content of Part 2.~in~Lemma~\ref{lem:delta and deltabar}.
We now recall that it was shown in~\cite{BloHofShp2017} that an arbitrary Phan amalgam over $\FF_q$ with $3$-spherical diagram possessing some non-trivial completion must satisfy property (D) and is isomorphic to $\amG(\delta)$ for some $\delta$, so the second claim follows from the first.
\epf

\medskip
As before, the Curtis-Tits amalgams fall into two categories: Orientable Curtis-Tits groups are essentially groups of Kac-Moody type (Theorem~\ref{thm:OCT realization}). Non-orientable Curtis-Tits groups can be obtained as (central extensions of) subgroups of groups of Kac-Moody type fixed under an involution that interchanges positive and negative roots groups and permutes types non-trivially (Theorem~\ref{thm:NOCT realization}). 
Note that the latter involutions are not Phan involutions. 
Completions of Phan amalgams are obtained as subgroups of groups of Kac-Moody type fixed under a Phan involution (which does fix the types) (Lemmas~\ref{lem:fixed CT amalgam is Phan amalgam}~and~\ref{lem:delta and deltabar}).
We show (Proposition~\ref{prop:G is a lattice of bhKM+}) that the completions of non-orientable Curtis-Tits groups are lattices in the ambient group of Kac-Moody type, thus generalizing a result from~\cite{GraMuh2008}.

Now suppose $\Gamma$ is non-spherical and non-affine and that $\fk$ is finite.
By~\cite{CapRem2009}, the corresponding orientable Curtis-Tits groups are almost simple (Corollary~\ref{cor:OCT almost simple}).
Note that the diagram of a non-orientable Curtis-Tits groups is either  $\tA_{n-1}$ or  it is non-spherical and non-affine.
In the present note we use a result from~\cite{CapHum2015} to show that non-orientable Curtis-Tits groups over finite fields with $3$-spherical diagram different from $\tA_{n-1}$ are acylindrically hyperbolic and in particular are not abstractly simple (Theorem~\ref{thm:non-orientable CT groups are not simple}).

In the case  where $\Gamma$ has type $\tA_{n-1}$ both orientable and non-orientable Curtis-Tits groups have interesting quotients~\cite{BloHof2014a}. In Subsection~\ref{subsec:examples} we also exhibit finite quotients of non-orientable Curtis-Tits groups with non-spherical and non-affine diagram.

Throughout the paper it is our intention to be as concrete and explicit as we can be. 

\section{Curtis-Tits and Phan amalgams}\label{sec:amalgams}

\subsection{Amalgams of Curtis-Tits and Phan type}\label{subsec:CTP amalgams}
%
%
\bde\label{dfn: CTP structure}
Let $\liediag=(I, E)$ be  a Lie  diagram.
A  {\em   Curtis-Tits (resp.~Phan)  amalgam with Lie diagram  $\liediag$ over $\FF_q$}  is a collection 
$\amG=\{\amgrpG_{i},\amgrpG_{i, j}, \amg_{i,j} \mid  i, j \in I\}$ 
 such that 
for every $i,j\in I$, $\amg_{i,j}\colon \amgrpG_i\to\amgrpG_{i,j}$ is a homomorphism of groups and, setting 
 $\bamgrpG_i=\amg_{i,j}(\amgrpG_i)$, the triple  
$(\amgrpG_{i,j}, \bamgrpG_i, \bamgrpG_j)$ is a Curtis-Tits / Phan  standard pair of type $\liediag_{i,j}(q^e)$, for some $e\ge 1$ as defined in~\cite{BloHofShp2017}.
Moreover $e=1$ is realized for some $i,j\in I$.
For any subset $K\sbe I$ , we let 
\begin{align*}
\amG_K&=\{\amgrpG_{i},\amgrpG_{i, j}, \amg_{i,j} \mid  i, j \in K\}. 
\end{align*}

A {\em completion} of $\amG$ is a group $\compG$ together with a collection  $\compg_\bullet=\{\compg_i,\compg_{i,j}\colon i,j\in I\}$ of homomorphisms $\compg_i\colon \amgrpG_i\to \compG$, and  $\compg_{i,j}\colon \amgrpG_{i,j}\to \compG$, whose images - often denoted ${\compG}_i=\compg_i(\amgrpG_i)$ - generate $\compG$, such that for any $i,j\in I$,  $\compg_{i,j}\after\amg_{i,j}=\compg_i$.
The amalgam $\amG$ is {\em non-collapsing} if it has a non-trivial completion.
As a convention, for any subgroup $\amgrp{H}\le \amgrpG_J$, let $\comp{H}=\compg(\amgrpH)\le \compG$.

A completion $(\ucompG,\ucompg_\bullet)$ is called {\em universal} if for any completion $(\compG,\compg_\bullet)$ there is a unique surjective group homomorphism $\umap\colon \ucompG\to \compG$ such that $\compg_\bullet=\umap\after\ucompg_\bullet$. A universal completion always exists and is unique, but it may be trivial.
\ede

\bde\label{dfn:standard amalgam}
Let $\liediag=(I, E)$ be  a Lie  diagram.
The   {\em standard  Curtis-Tits (resp.~Phan)  amalgam with Lie diagram  $\liediag$ over $\FF_q$}  is the Curtis-Tits (resp.~Phan) amalgam 
$\famG=\{\amgrpG_{i},\amgrpG_{i, j}, \famg_{i,j} \mid  i, j \in I\}$  in which $\famg_{i,j}$ is the standard identification map as defined in~\cite{BloHofShp2017} for all $i,j\in I$.
\ede

We now describe all Curtis-Tits and Phan amalgams arising from the classification results of~\cite{BloHofShp2017}.
To this end, first consider certain groups $\amgrpC_i$ of automorphisms of the vertex groups $\SL_2(q)$ (Curtis-Tits case) and $\SU_2(q)$ (Phan case). These are certain subgroups of the vertex groups of the Coefficient system $\amA$ of~\cite{BloHofShp2017}.
\paragraph{$\amgrpC_i$ Curtis-Tits case} 
For any $i\in I$, there is some $e\in \NN$ such that $\amgrpG_i=\SL_2(q^e)$ via the standard identification maps and $\amgrpC_i=Aut(\FF_{q^e})\times\langle \tau\rangle$ (with $\tau$ of order $2$).
Here $\alpha\in \Aut(\FF_{q^e})$ acts as a Frobenius automorphism and $\trin$ acts as transpose-inverse.

\paragraph{$\amgrpC_i$ Phan case} 
For any $i\in I$, we have  $\amgrpG_i=\SU_2(q)$ via the standard identification maps and $\amgrpC_i=\Aut(\FF_{q^2})$.
Here $\alpha\in \Aut(\FF_{q^2})$ acts as a Frobenius automorphism with respect to an orthonormal basis for the hermitian form.
Note that  $\trin$ (transpose-inverse) acts as $\sigma\colon x\mapsto x^q$ ($x\in \FF_{q^2}$) on $\SU_2(q)$.

\paragraph{The spanning tree for $\liediag$}
As customary we view the Dynkin diagram $\liediag$ as an oriented edge-labeled graph. Let $\uliediag$ denote the underlying undirected simple graph.
We now fix a spanning tree  $T=(I,E(T))$ for $\uliediag$ and let $\edg\uliediag-\edg T=\{\{i_s,j_s\}\colon s=1,2,\ldots,r\}$ together with certain integers $\{e_s\colon s=1,2,\ldots,r\}$.
In the Phan case, any spanning tree suffices and $e_s=1$ for all $s$. In the Curtis-Tits case select $T$ such that (see~\cite{BloHofShp2017}):
 \begin{enumerate}
 \item\label{cond:A2} $(\amgrpG_{\{i_s,j_s\}},\amg_{i_s,j_s}(\amgrpG_{i_s}),\amg_{i_s,j_s}(\amgrpG_{j_s}))$ has type $A_2(q^{e_s})$, where $e_s$ is some power of $2$.
  \item\label{cond:minimal e} There is a loop $\Lambda_s$ containing $\{i_s,j_s\}$ such that 
  any vertex group of $\Lambda_s$ is isomorphic to $\SL_2(q^{e_s 2^l})$ for some $l\ge 0$.  
 \end{enumerate}
\bde\label{dfn:representative amalgams}
 The main results of~\cite{BloHofShp2017} now says that the Curtis-Tits (resp.~Phan) amalgams over $\FF_q$ with diagram $\liediag$ are, up to type preserving isomorphism, in bijection with the set 
  \begin{align*}
  \amgrpC=\prod_{s=1}^r\amgrpC_{i_s}.
  \end{align*}
Under this bijection, the sequence $\delta=(\delta_s)_{s=1}^r$ corresponds to the amalgam 
\begin{align}
\amG(\delta)&=\{\amgrpG_{i},\amgrpG_{i, j}, \amg^\delta_{i,j} \mid  i, j \in I\} \label{eqn:amG^delta}
\end{align} 
where $\amg^\delta_{i,j}=\famg_{i,j}$ for all $i,j\in I$ except that 
$\amg_{i_s,j_s}^\delta= \famg_{i_s,j_s}\after\delta_s$ for all $s=1,2,\ldots,r$.
We shall call $\famG(\delta)$ {\dfn  the amalgam representing $\delta$}.

We finish this discussion with some terminology applying only to Curtis-Tits amalgams. We say that $\famG(\delta)$ is {\dfn orientable} if  $\delta_s\in \Aut(\FF_q^{e_s})$ for all $s=1,\ldots,r$
 (that is, $\delta$ does not involve any $\trin$) and {\dfn non-orientable} otherwise.
Note that we can interpret $\delta$ as the image of a homomorphism
\begin{align*}
\classmap\colon \pi_1(\liediag,0)& \to \amgrpC\\
[\Lambda_s]&\mapsto \delta_s
\end{align*}
where $0$ is some base vertex of $\liediag$ and $[\Lambda_s]$ denotes the homotopy class of the loop $\Lambda_s$ above.
Consider the composition of $\classmap$ and the natural projection map:
\begin{align}
\classmap^*\colon \pi_1(\liediag,0)\to \langle \trin\rangle\cong \ZZ/2\ZZ\label{eqn:classmap star}
\end{align}
Then $\amG$ is orientable if and only if the image of the corresponding $\classmap^*$ is trivial.
\ede
From now on $\amG=\famG(\delta)=\{\amgrpG_{i},\amgrpG_{i, j}, \amg_{i,j} \mid  i, j \in I\}$ 
 for some $\delta$.
We shall only consider Curtis-Tits amalgams $\amG$ over a finite field $\FF_q$ possessing a weak system of fundamental root groups with $3$-spherical (but not spherical) diagram having no subdiagrams of type $C_2(2)$. 
The latter condition was not necessary for classification, but it is necessary when considering completions, as it is necessary to satisfy condition (co) of~\cite{Mu1999}.

\section{Orientable Curtis-Tits groups are groups of Kac-Moody type}\label{sec:curtis-tits realizations}
\subsection{Realization of Orientable Curtis-Tits amalgams} Realization of all Curtis-Tits amalgams arising from the classification can be achieved along the lines of~\cite{BloHof2016,BloHof2014b}, where they were obtained for Curtis-Tits amalgams over a finite field $\FF_q$  with $q\ge 4$ satisfying property (D)  and having $3$-spherical simply-laced diagram. The present situation only requires us to modify the proof in certain places, and we will content ourselves with pointing out these differences.

We start by assuming that $\amL=\{\amgrpL_i,\amgrpL_{i,j},\aml_{i,j}\colon i\in I\}$ is an orientable Curtis-Tits amalgam over $\FF_q$ with $3$-spherical diagram $\covdiag$ without $C_2(2)$-subdiagrams..
Following Tits~\cite{Ti1992} a {{\em group of Kac-Moody type}}{} is by definition a group with RGD  such that a central quotient is the subgroup of $\Aut(\De)$ generated by the root groups of an apartment in a Moufang twin-building {$\De$}{}. This central quotient will be called the associated {\em adjoint} group of Kac-Moody type. In this section we shall prove the following.
As a general reference for groups with root group datum we will use~\cite{CapRem2009a}. In particular, we assume that such a groups is generated by the root groups of the root group datum.

\bth \label{thm:OCT realization}
The universal completion of $\amL$ is a group of Kac-Moody type (and $\amL$ is the Curtis-Tits amalgam for this group) if and only if $\amL$ is orientable.
\eth

Theorem~\ref{thm:OCT realization} generalizes Corollary 1.2 of~\cite{BloHof2014b} and its proof follows the steps detailed in Section 5 of~\cite{BloHof2014b}.
So as not to repeat that proof nearly verbatim, we merely indicate how the more general assumptions of Theorem~\ref{thm:OCT realization} still yield the same result.

\medskip
\bpf(of Theorem~\ref{thm:OCT realization})
Following Subsection~5.1 of~\cite{BloHof2014b} we consider a simply-connected locally split Kac-Moody group over $\FF_q$ with diagram $\covdiag$ and consider the twin-building $\De=((\De_+,\delta_+),(\De_-,\delta_-),\delta_*)$ associated to its twin BN-pair $(B^+,N,B^-)$. Now as $\covdiag$ is $3$-spherical and has no $C_2(2)$ subdiagrams, it satisfies condition (co) of~\cite{MuRo1995} so the local structure determines the global structure. For the same reason $\covdiag$ satisfies condition (${\rm co}^*$) from~\cite{Cap2007} so that by the main result of that paper, we obtain the group as a central quotient of the universal completion of an amalgam
\begin{align*}
\cR&=\{R_i,R_{i,j}, \rho_{i,j}\mid i,j\in I\},
\end{align*}
where $R_i=\langle X_i^+,X_i^-\rangle$, $R_{i,j}=\langle R_i,R_j\rangle$,  
and $\{X_i^+\mid i\in I\}$ is a selection of positive root groups corresponding to a fundamental system of positive roots, and $ \rho_{i,j}\colon R_i \into R_{i,j}$ is given by inclusion of subgroups in $G$.
Now  $\cR$ is the desired Curtis-Tits amalgam. Note that the analog of Lemma 5.2 is not valid as for instance the center of $\Sp_4(q)$ meets one of the rank $1$ Levi groups, but this does not affect the validity of the conclusion.
Note that $\Chi=\{X_i^+,X_i^-\colon i\in I\}$ is the weak system of fundamental root groups and the fact that we can select the $+$ signs to correspond to a fundamental system of positive roots means that the amalgam is orientable.

We now show that the universal completion of any orientable $\amL$ has a central quotient that is the automorphism group of a twin-building following Subsection 5.3 of~\cite{BloHof2014b}. We now take the definition of a sound Moufang foundation as in~\cite{Mu1999}. As $\amL$ is orientable, we obtain the rank-$2$ Moufang buildings $\De_{i,j}$ of type $\covdiag_{i,j}$ of the foundation $\sF$ from the Curtis-Tits standard pairs $(\amgrpL_{i,j},\amgrpL_i,\amgrpL_j)$ in $\amL$ using the Borel subgroups $\{\amgrpB_{i,j}^+,\amgrpB_{i,j}^-\}$, which are uniquely determined by $\{\amgrpX_i^\pm,\amgrpX_j^\pm\}\sbe \Chi$. We now select the chamber $C_{i,j}\in \De_{i,j}$ of the foundation to be those associated to $\amgrpB_{i,j}^+$. Also for each $i\in I$ we  select an auxiliary chamber $C_i$ in the rank-$1$ building $\De_i$ associated to the BN pair $(\amgrpB_i^+,N_i,\amgrpB_i^-)$. We define an inclusion map $\theta_i^j\colon \De_i\to \De_{i,j}$ induced by $\aml_{i,j}$ and let the restriction maps of $\sF$ be given by $\theta_{i,j}^{j,k}=\theta_j^k\after(\theta_j^i)^{-1}$.
It is immediate that the $\theta_{i,j}^{j,k}$ satisfy the condition (Fo3) of a Moufang foundation.
Soundness of $\sF$ follows from the fact that $\covdiag$ is $3$-spherical and the fact that the subamalgam of $\amL$ of type $\covdiag_J$ is the Curtis-Tits amalgam associated to the corresponding spherical building.
The fact that the signature of $\Chi$ determines a twin-apartment is proved as in~\cite{BloHof2014b}. 
The proof is completed exactly as in~\cite{BloHof2014b}: By~\cite{Mu1999}, the sound Moufang foundation can be  integrated to a twin-building $\De$ whose automorphism group contains a non-trivial homomorphic image of the original Curtis-Tits amalgam, generated by the root groups associated to the roots intersecting the $E_2(c)$ for some chamber $c$. Universality gives a homomorphism from the universal completion $\ucompL$ of $\amL$ to the subgroup $\Aut(\De)^{\dagger}$ of $\Aut(\De)$ generated by these root groups; the kernel of the action of $\ucompL$ on $\De$ must be central, as required.
\epf

\bre
In~\cite{BloHof2016, BloHof2014a} we obtained rather precise information on the particular central subgroups appearing in the kernel of the map from the universal completion of $\amL$ to $\Aut(\De)$. For the purposes of the present paper, it suffices to establish the existence of a completion for $\amL$.
However, we expect that the techniques developed in~\cite{BloHof2016,BloHof2014a} can be used to handle this more general case, although the details will probably a bit more involved.
\ere

In subsequent sections we shall prove that non-orientable Curtis-Tits (resp.~Phan) amalgams have non-trivial completions inside the subgroup of an orientable Curtis-Tits group  fixed under a certain Cartan (resp.~Phan) involution. 
\subsection{The twin-building $\De$ associated to $\ucompL$}\label{subsec:twin-building of amL}
Let $(\ucompL,\ucompl)$ be the universal completion of $\amL$.
For future reference, we introduce the notation necessary to talk about the twin-building related to $\ucompL$.
Note that $\ucompL$ is a group of Kac-Moody type over $\FF_q$ with diagram $\covdiag$.
Let $(W,S=\{s_i\colon i\in I\})$ the associated Coxeter system with root system $\Phi$.
Now $\ucompL$ is a group with a locally finite root group datum $\{\compU_\alpha\colon \alpha\in \Phi\}$ (namely $\compU_\alpha$ is finite for all $\alpha\in \Phi$ (see~\cite{CapRem2009,Cap2007}). 
This means in particular that  $\ucompL=\langle \compU_\alpha\colon \alpha\in \Pi\rangle$ (for a root base $\Pi$ of $\Phi$) has a twin $BN$-pair  $((\compB^+,\compN), (\compB^-,\compN))$, where $\compB^\vep=\compD\compU^\vep$, setting 
\begin{align*}
\compU^\vep&=\langle \compU_\alpha\colon \alpha\in \Phi^\vep\rangle &&  (\vep=+,-) \\
\compN&=\langle \mu(u)\colon u\in \compU_\alpha-\{1\}, \alpha\in \Pi\rangle,\\
\compD&=\bigcap_{\alpha\in \Phi}N_{\ucompL}(\compU_\alpha).
\end{align*}
In fact this is the twin BN-pair giving rise to a twin-building $((\De_+,\delta_+),$ $ (\De_-,\delta_-),\delta_*)$.
As proved in~\cite{BloHofShp2017} the weak system of fundamental root groups can be selected so that for some  fundamental system $\Pi=\{\alpha_i\colon i\in I\}$ of $\Phi$, and all $i\in I$, we have 
$\compU_{\alpha_i}=\compl(\amgrpX_i^+)$.

\section{Non-orientable Curtis-Tits groups are fixed groups of Cartan involutions}\label{sec:NO CT groups}

\subsection{The ambient orientable amalgam $\amL$}
We shall now assume that $\amG=\famG(\delta)=\{\amgrpG_{i},\amgrpG_{i, j}, \amg_{i,j} \mid  i, j \in I\}$ 
 where $\delta=(\delta_s)_{s=1}^r\in \prod_{s=1}^r \amgrpC_{j_s}$ is non-orientable. This means that the map $\classmap^*$ of~\eqref{eqn:classmap star} is surjective.
\bde {(The covering of diagrams $p\colon \covdiag\to\liediag$)}
Consider the map $\classmap^*\colon \Pig\to \langle \tau\rangle$.
Its kernel is the fundamental group of a two-sheeted covering $p\colon \covdiag\to \liediag$ sending some vertex $\hat{0}$ to $0$.
Since $\liediag$ has no circuits of length $\le 3$, the quotient $\pi(\Gm,0)/\pi(\covdiag,\hat{0})=\langle\theta\rangle\cong \ZZ/2\ZZ$ acts as a group of deck transformations 
 commuting with $p$; in particular $p$ does not fix points or edges.
%
%
\ede

We now lift $\amG$ to a locally isomorphic amalgam $\amL$ defined over $\covdiag$ and extend $\theta$ to $\amL$.
\bde\label{dfn:hamG}
Let $\amL=\{\bL_{i},\bL_{i, j}, \hpsi_{i,j} \mid  i, j \in \hI\}$ be the amalgam such that, for all $i,j\in \hI$, 

\begin{enumerate}
\item[($\amL$1)] 
$\bL_i\mbox{  is a copy of }\G_{p(i)}$,
\item[($\amL$2)] 
$ \bL_{i,j}\mbox{  is a copy of }
\begin{cases}
\G_{p(i),p(j)} &\mbox{if  } \{i,j\}\in \hE, \\
\G_{p(i)}\times \G_{p(j)} & \mbox{else. } \\
\end{cases}$
\item[($\amL$3)] 
$ \hpsi_{i,j}
 =
 \begin{cases}
\vph_{p(i),p(j)} & \mbox{if } \{i,j\}\in \hE,\\
 \mbox{ canonical inclusion } & \mbox{else.}
\end{cases}$
\end{enumerate}

This means the following. Fix some $J\sbe \hI$ with $1\le |J|\le 2$, and denote by $\pi\colon \amL\to \amG$ the homomorphism of amalgams induced by $p$.
That is, identifying $\G_{p(J)}$ with its copies $\bL_J$ and $\bL_{\theta(J)}$, we let the maps $\pi_J$ and $\pi_{\theta(J)}$ be the identity mappings.
Then, we have a commuting diagram of isomorphisms

\begin{align}\label{eqn:hamG covers amG}
\xymatrix{
x\in \bL_J\ar[dr]^{\pi_J} \ar[rr]^\theta && \bL_{\theta(J)} \ni x \ar[dl]^{\pi_{\theta(J)}}\\
                            & x\in \G_{p(J)} & 
}
\end{align}
Thus, in~\eqref{eqn:hamG covers amG}, also $\theta$ is given by the identity mapping.
Then, by ($\amL$3) $\theta$ is an automorphism of $\amL$. Also $\pi$ can be viewed as a $2$-covering of Curtis-Tits amalgams.
\ede
Since $\pi$ induces an isomorphism on every vertex and edge group of $\amL$ we have the following.
If $\delta_s\in\Aut(\FF_{q^{2^{e_s}}})$, then the $p$-fiber over $\La_s$ consists of two disjoint loops, and the subamalgams of $\amL$ induced on each of these is isomorphic to the one induced on $\La_s$, hence they correspond to $\delta_s$ as well.
Otherwise the $p$-fiber over $\La_s$ is a single loop $\hat{\La}_s$ doubly covering $\La_s$ and it corresponds to $\delta_s^2\in\Aut(\FF_{q^{2^{e_s}}})$.
Hence $\amL$ is orientable. It follows from the classification theorem that it is isomorphic to some standard orientable Curtis-Tits amalgam.

\subsection{Realization of $\amG$ in the twisted group of Kac-Moody type $\amgrpL^\theta$}
Since $\amL$ is orientable, Theorem~\ref{thm:OCT realization} implies that $\amL$ has a completion $(\compL,\compl)$, where $\compL$ is a group of Kac-Moody type. Hence also the universal completion $(\ucompL,\ucompl)$ of $\amL$ is not trivial.
For convenience, we shall replace $\amL$ by its image in $\ucompL$. 
Note that this image does not have to be isomorphic to $\amL$, but when considering completions, there is no loss of generality.

By universality, the automorphism $\theta\colon \amL\to\amL$ induces a Cartan involution, also denoted $\theta$, on the universal completion $\ucompL$.  
\bde\label{dfn:Cartan involution}
A {\dfn Cartan involution} of a group with twin- $BN$-pair $((\compB^+,\compN),(\compB^-,\compN))$ is an automorphism $\theta$ satisfying   
\begin{enumerate}
\item $\theta^2=\id$
\item$(\compB^+)^\theta=\compB^-$,
\item $\theta$ normalizes the Weyl group $W=\compN/\compB^-=\compN/\compB^-$ inducing a graph automorphism of $\liediag$ without fixed vertices or edges.
\end{enumerate}
\ede
\noindent
Let $\ucompL^\theta$ be the fixed group of $\ucompL$ under $\theta$ and define the amalgam of fixed subgroups
\begin{align*}
\amL^\theta & = \{\aml^\theta_{i,j}\colon \langle \amgrpL_i,\amgrpL_{\theta(i)}\rangle^\theta\into\langle\amgrpL_{i,j},  \amgrpL_{\theta(i),\theta(j)}\rangle^\theta\colon i,j\in I\}, 
\end{align*}
where, for each $i,j\in I$, we consider the fixed subgroups  
\begin{align*}
\langle \amgrpL_i,\amgrpL_{\theta(i)}\rangle^\theta &=\{x x^\theta\colon x\in \amgrpL_i\},\\
\langle \amgrpL_{i,j},\amgrpL_{\theta(i), \theta(j)}\rangle^\theta &=\{x x^\theta\colon x\in \amgrpL_{i,j}\}.
\end{align*}
There is a non-trivial surjective morphism $\amG\to \amL^\theta$.
In~\cite{BloHof2016} this is shown for the image of $\amL^\theta$  in $\compL$. 
However, the proof translates verbatim to obtain the result in the present more general setting.
Hence there is a non-trivial completion $(\compG,\compg)$ of $\amG$ with $\compG=\langle\amL^\theta\rangle \le \ucompL^\theta$. 
Thus, we have obtained the following.
\bth\label{thm:NOCT realization}
Let $\amG$ be a non-orientable Curtis-Tits amalgam over $\FF_q$ with connected $3$-spherical Dynkin diagram $\liediag$ having no subdiagrams of type $C_2(2)$. Then, there is a group of Kac-Moody type $\ucompL$ over $\FF_q$ whose diagram is a two-sheeted covering of $\liediag$ equipped with a Cartan involution $\theta$ such that the fixed group $\ucompL^\theta$ contains a non-trivial completion of $\amG$.
\eth

\bre
Let $(\ucompG,\ucompg)$ be the universal completion of $\amG$.
Then, there is a surjective homomorphism
\begin{align*}
\tilde{\compG}\onto \compG\le \ucompL^\theta\le \ucompL.
\end{align*}
Let $\tilde{\pi}\colon\ucompL\to\compL$ be the map of completions of $\amL$ given by universality. 
Now $\ucompL$ is a central extension of $\compL$. If in fact $\compL=\ucompL/Z(\ucompL)$, then $\theta$ induces an involution of $\compL$ so that $\tilde{\pi}$ sends $\ucompL^\theta$ to $\compL^\theta$. In this terminology,~\cite{BloHof2016} investigates the index $[\compL^\theta\colon\tilde{\pi}(\compG)]$ in great detail in the case of simply-laced diagrams.
\ere
\section{Realization of Phan amalgams}\label{sec:Realization of Phan amalgams}
\subsection{Phan involutions}\label{subsec:phan involutions}
Recall the following definition of a Phan involution for groups with twin-root group datum.
\bde\label{dfn:orientable phan involution}
A {\dfn Phan involution} $\theta$ of a group $\compL$ with a twin root group datum (see Subsection~\ref{subsec:twin-building of amL}) is an automorphism of $\compL$ such that 
\renewcommand{\theenumi}{\roman{enumi}}\begin{enumerate}
\item\label{Phan Inv i} $\theta^2=\id$
\item\label{Phan Inv ii} $(\compB^+)^\theta=\compB^-$,
\item\label{Phan Inv iii} $\theta$ centralizes the Weyl group $W=\compN/\compD$.
\end{enumerate}
\ede

Now let $\amG=\{\amg_{i,j}\colon\amgrpG_i\to  \amgrpG_{i,j}\colon i,j\in I\}$ be a (possibly non-orientable) Curtis-Tits amalgam over $\FF_q$ with $3$-spherical diagram $\liediag$.
Let  $(\ucompG,\ucompg)$ denote its universal completion and let $(\compG,\compg)$ be some completion 
 with canonical map $\uammapp\colon \ucompG\to\compG$.
 
\bde\label{dfn:phan involution of CT group}
A {\dfn Phan involution}  of $\amG$ is an automorphism $\theta=\{\theta_i,\theta_{i,j}\colon i,j\in I\}$ of $\amG$ that induces a Phan involution on each group of $\amG$.
\ede

From now on let $\theta=\{\theta_i,\theta_{ij}\colon i,j\in I\}$ be a Phan involution of the Curtis-Tits amalgam $\amG$.
As a source of examples, we have the following observation.
\ble\label{lem:phan involutions inducing phan involutions on CT amalgam}
Suppose that $\amL$ is the Curtis-Tits amalgam arising from the action of a group of Kac-Moody type $\compL$ on its twin-building $\De$ over $\FF_q$ with $3$-spherical diagram $\covdiag$ and that 
 $\compl\colon \amL\to\compL$ is the completion map. 
 If $\theta$ is a Phan involution of $\compL$ preserving
 $\amL$, then it induces a Phan involution on $\amL$.
 Conversely, a Phan involution $\theta$ of $\amL$ induces a Phan involution on $\ucompL$.
\ele

\bpf
First note that any automorphism of $\amL$ induces a unique automorphism of its universal completion $(\ucompL,\ucompl)$ that preserves $\amL$. 
Conversely, any automorphism of $\compL$ which preserves $\amL$,  induces an automorphism of $\amL$.
In both cases, if the original automorphism is an involution, then so is the induced one.
Let $\theta$ denote the automorphism of $\amL$ as well as the induced automorphism of $\compL$ in the forward direction, and of $\ucompL$ in the backward direction.

As we saw in~Subsection~\ref{subsec:twin-building of amL} not only $\compL$, but also $\ucompL$ is a group with root group datum. So as not to overload notation, in both cases (directions) we shall denote the datum
 $\{\compU_\alpha\colon \alpha\in \Phi\}$. 

Note that when proving either direction, we can view $\amL$ as a concrete amalgam for $\compL$  (resp.~$\ucompL$) in the sense that for any $i,j\in I$, the connecting map 
 $\aml_{i,j}$ (resp.~$\uaml_{i,j}$)  is just inclusion of subgroups 
(In the backward direction, for purposes of completions, it is harmless to identify $\amL$ with its image in $\ucompL$). In the terminology of~\cite[Section 4]{BloHofShp2017} this means that 
  if $\Chi=\{\amgrpX_i^+,\amgrpX_i^-\colon i\in I\}$ is the weak system of fundamental root groups of $\amL$, then 
for any $j\in I$, we have $\amgrpX_i^\vep=\aml_{i,j}(\amgrpX_i^\vep)$ (resp.~ $\amgrpX_i^\vep=\uaml_{i,j}(\amgrpX_i^\vep)$) ($\vep=+,-$). By uniqueness of fundamental root groups in $\amgrpL_{i,j}$ whenever $\covdiag_{i,j}\ne A_1\times A_1$ and connectedness of $\covdiag$, we then have $\Chi=\{\compU_{\alpha_i},\compU_{-\alpha_i}\colon i\in I\}$.  The significance here is that $\Pi=\{\alpha_i\colon i\in I\}$ is a root basis for $\Phi$, so that 
\begin{align}
W\Pi=\Phi.\label{eqn:Chi^+ is a root basis for Phi}
\end{align}
Since in the forward (resp.~backward) direction the automorphism $\theta$ of $\compL$ (resp.~of~$\ucompL$) preserves each group of $\amL$,  it follows that $\theta$ acts as 
\begin{align}
\compU_{\alpha}^\theta&=\compU_{-\alpha}\mbox{ for all }\alpha\in \Pi \label{eqn:theta swaps fundamental root groups}
\end{align}

{We now complete the proof of the forward direction. 
Combining property~\eqref{Phan Inv iii} of the Phan involution $\theta$ on $\compL$, with~\eqref{eqn:Chi^+ is a root basis for Phi},~\eqref{eqn:theta swaps fundamental root groups} and applying the RGD axioms, we see that $\theta$ preserves the root group datum, while interchanging $\{\compU_\beta\colon \beta\in \Phi_J^+\}$ and  $\{\compU_\beta\colon \beta\in \Phi^-_J\}$ for any $J\sbe I$; in particular, $\theta$ preserves $\compN_J$ and $\compD$, 
 and centralizes $W_J=\compN_J/\compD$.
Hence, $\theta_i$ and $\theta_{i,j}$ satisfy properties~\eqref{Phan Inv ii}~and~\eqref{Phan Inv iii} for each $i,j\in I$.

We now establish the backward direction in a similar manner using~\eqref{eqn:Chi^+ is a root basis for Phi}~and~\eqref{eqn:theta swaps fundamental root groups}. Thus, in order to establish that the automorphism $\theta$ of $\ucompL$ has property~\eqref{Phan Inv ii} it suffices to show that 
it has property~\eqref{Phan Inv iii}. Recall that $N$ is generated by elements $\mu(u_i)$, for some $u_i\in \compU_{\alpha_i}$ with $\alpha_i\in \Pi$. 
Let $(W,\{s_i\colon i\in I\})$ denote the Coxeter system where $s_i=\mu(u_i)D$.
Now note that $(\ucompL_i,\{\compU_{\alpha_i},\compU_{-\alpha_i}\})$ is a group with root group datum since $\ucompL_i=\langle \compU_{\alpha_i},\compU_{-\alpha_i}\rangle_{\ucompL}$ (cf.~\cite[\S2.3]{CapRem2009a}). 
Since $\theta_i$ satisfies property~\eqref{Phan Inv iii}, the element $\mu(u_i)^{\theta_i}\in \compU_{\alpha_i}u_i'\compU_{\alpha_i}$ (with $u_i'\in \compU_{-\alpha_i}$) must conjugate $\compU_\beta$ to $\compU_{s_i\beta}$ for all $\beta\in \Phi$; in particular it must do so for $\beta=\pm\alpha_i$.
However, in this standard root group datum of $\SL_2(\fk)$ for some field $\fk$, one verifies that this means that 
 $\mu(u_i)^{\theta_i}=\mu(u_i')$. Clearly $\mu(u_i)D=s_i=\mu_(u_i')D$. Thus $\theta$ centralizes $W$, as required.
\epf
}

\subsection{Fixed subamalgams}\label{subsec:fixed subamalgams}

\bde\label{dfn:fixed amalgam}
Let $\theta$ be a Phan involution of the Curtis-Tits amalgam $\amG$. Define the fixed amalgam of $\amG$ under $\theta$ as $\amF=\{\amgrpF_{i,j},\amgrpF_i, \amf_{i,j}\colon i,j\in I\}$, 
 where 
 \begin{align*}
 \amgrpF_i & = \amgrpG_i^{\theta_i}, \\
 \amgrpF_{i,j} & = \amgrpG_{i,j}^{\theta_{i,j}},\\
 \amf_{i,j} & = \amg_{i,j}|_{\amgrpF_i}.
 \end{align*}
We may therefore unambiguously write $\amF=\amG^\theta$, $\amgrpF_J=\amgrpG_J^\theta$ ($J\sbe I$ with $0<|J|\le 2$) and denote the inclusion maps simply with $\amg_{i,j}$ instead of $\amf_{i,j}$.
\ede
\bre\label{rem:orthogonal amalgams}
The case where $\theta$ induces $\trin$ on all groups of $\amG=\famG$ is the one studied in~\cite{CapHum2015}. 
For $A_n$ diagrams the situation is studied in detail in~\cite{Hoffman:2013aa}.
In the present paper we are mostly interested in the case where $\theta$ induces $\trin$ together with some field involution.
\ere
\medskip
\noindent

%
\medskip
We now establish the relation between $\amG$ and $\amF$ in terms of the representative amalgams of Definition~\ref{dfn:representative amalgams}.
Let $\famG$ be the standard Curtis-Tits amalgam over $\FF_{q^2}$ with connected $3$-spherical diagram $\liediag$ (as in Definition~\ref{dfn:standard amalgam}), and assume that that for any $i,j\in I$,
$\liediag_{i,j}\in \{A_1\times A_1, A_2,C_2/B_2\}$.
Note that since $\liediag$ has no subdiagrams of type $\twA_3$, in fact we have that, for any $i,j\in I$
$(\amgrpG_{i,j},\bamgrpG_i,\bamgrpG_j)$ is a Curtis-Tits standard pair of type $\liediag_{i,j}(q^2)\in \{A_1(q^2)\times A_1(q^2), A_2(q^2),C_2(q^2)/B_2(q^2)\}$.

Fix some spanning tree $\Sigma\sbe \liediag$ and suppose that $\edg\liediag-\edg\Sigma=\{\{i_s,j_s\}\colon s=1,2,\ldots,r\}$ so that $H_1(\covdiag,\ZZ)\cong\ZZ^r$. 
Let $\amG=\famG(\delta)$ for some $\delta\in \prod_{s=1}^r\amgrpC_{i_s}$.

Note that $\amgrpC_{i_s}=\Aut(\FF_{q^2})\times\langle\trin\rangle$ for all $s$.
We let $\theta=\sigma\after\trin$, where we recall that $\sigma\colon x\mapsto x^q$ for $x\in \FF_{q^2}$.

\ble\label{lem:fixed CT amalgam is Phan amalgam}
The fixed amalgam $\amF$ of $\amG$ under $\theta$ is a Phan amalgam over $\FF_q$ with diagram $\liediag$.
\ele

\bpf
We simply have to verify that if $(\amgrpG_{i,j},\bamgrpG_i,\bamgrpG_j)$ is a Curtis-Tits standard pair of type $\liediag_{i,j}(q^2)\in \{A_1(q^2)\times A_1(q^2), A_2(q^2),C_2(q^2)/B_2(q^2)\}$, then 
$(\amgrpF_{i,j}^{\theta_{ij}},\bamgrpF_i^{\theta_{i}},\bamgrpF_j^{\theta_{j}})$ is a Phan standard pair of type $\liediag_{i,j}(q)\in \{A_1(q)\times A_1(q), A_2(q),C_2(q)/B_2(q)\}$.
Now note that $\theta$ is defined with respect to the standard basis used to define the Curtis-Tits standard pair.
Thus, the fixed group under $\theta$ is the intersection of the group in the Curtis-Tits standard pair and the unitary group preserving the hermitian form for which this standard basis is orthonormal.
In every case this is exactly the corresponding vertex or edge group of the Phan standard pair with the same diagram as the Curtis-Tits standard pair.
\epf

Let $\famF$ be the standard Phan amalgam with diagram $\liediag$ over $\FF_q$. Let $\amF$ be the fixed amalgam of $\amG$ under $\theta$.  By Lemma~\ref{lem:fixed CT amalgam is Phan amalgam} $\amF$ is a Phan amalgam. By the classification this means that $\amF=\famF(\bar{\delta})$ for some $\bar{\delta}$. We now identify $\bar{\delta}$ given $\delta$.
\ble\label{lem:delta and deltabar}
\begin{enumerate}
\item Suppose $\amG=\famG(\delta)$ and $\amF=\famF({\bar{\delta}})$ Then,   
$\bar{\delta}$ is the image of $\delta=(\delta_{s})_{s=1}^r$ under the map
\begin{align}
\prod_{s=1}^r \amgrpC_{i_s} & \to \prod_{s=1}^r \amgrpC_{i_s}/\langle \theta_{i_s}\rangle.\label{eqn:classification restriction}
\end{align}
\item As a consequence, every Phan amalgam $\amF$ appears as the fixed amalgam in $2^r$ pairwise non-isomorphic Curtis-Tits amalgams $\amG$.
\end{enumerate}
\ele
\bpf
Note that in the case of Curtis-Tits as well as Phan standard pairs, the standard identification maps, i.e.~the connecting maps for $\famG$ and $\famF$ are "identity maps''. 
It follows that 
\begin{align}
\famG^\theta&=\famF.\label{eqn:famGtheta=famF}
\end{align}
To conclude 1.~from~\eqref{eqn:famGtheta=famF} it suffices to note that if $\famF(\bar{\delta})$ is the fixed amalgam of $\famG(\delta)$, then 
\begin{align*}
\famg_{i_s,j_s}|_{\amgrpG_{i_s}^\theta}\after \bar{\delta}_{i_s}=\famf_{i_s,j_s}\after  \bar{\delta}_{i_s}=\amf_{i_s,j_s} & = \amg_{i_s,j_s}|_{\amgrpG_{i_s}^\theta}=\famg_{i_s,j_s}|_{\amgrpG_{i_s}^\theta}\after \delta_{i_s}|_{\amgrpG_{i_s}^\theta},
\end{align*}
 and by restricting to ${\amgrpG_{i_s}^\theta}$, we pass from 
$\amgrpC_{i_s}=\Aut(\FF_{q^2})\times \langle \trin\rangle$ to $\amgrpC_{i_s}/\langle\theta\rangle=\Aut(\FF_{q^2})$ by identifying $\trin=\sigma$.

Part~2.~follows from the classification of Curtis-Tits and Phan amalgams combined with the observation that 
 the map~\eqref{eqn:classification restriction} has a kernel of order $2^r$.
\epf

\subsection{Completions of Phan amalgams}\label{subsec:completions of Phan amalgams}
Note that, if $(\ucompG,\ucompg)$ is non-trivial, then so is the image of $\amF$ in $\ucompG$ under $\ucompg$.
In particular, $\amF$ also has a non-trivial completion contained in the fixed group $\ucompG^\theta$.
This proves the following:

\bpr\label{prop:Phan realization}
Let $\amF$ be a Phan amalgam over $\FF_q$ with connected $3$-spherical Dynkin diagram $\liediag$ having no subdiagrams of type $C_2(2)$.  Then, there is a  Curtis-Tits amalgam  $\amG$ over $\FF_{q^2}$ with diagram $\liediag$  whose universal completion $(\ucompG,\ucompg)$ equipped with a 
 Phan involution $\theta$ such that its fixed group $\ucompG^\theta$ contains a non-trivial completion $(\compF,\compf)$ of $\amF$.
\epr

\bco\label{cor:Phan realization in KM group}
Let $\amF$ be a Phan amalgam over $\FF_q$ with connected $3$-spherical Dynkin diagram $\liediag$ having no subdiagrams of type $C_2(2)$. 
Then, there is a group of Kac-Moody type $\ucompL$ over $\FF_{q^2}$ with diagram 
$\liediag$ equipped with a Phan involution $\theta$ such that the fixed group $\ucompL^\theta$ contains a non-trivial completion $(\compF,\compf)$  of $\amF$.
\eco

We note that the completion $(\compF,\compf)$ in Proposition~\ref{prop:Phan realization}~and~Corollary~\ref{cor:Phan realization in KM group} can be a proper subgroup of $\ucompG^\theta$ and the completion does not have to be universal.
\section{Non-orientable Curtis-Tits groups are lattices}
\subsection{Cartan involutions of the twin-building $\De$}
We continue the notation from Section~\ref{sec:NO CT groups}.
Note that the properties of a Cartan involution can be reformulated in terms of the action on the building, as follows.
\ble\label{lem:theta on Delta}
Let $\theta$ be a Cartan involution as in Definition~\ref{dfn:Cartan involution}. Then $\theta$ induces an automorphism, also denoted $\theta$, on $\De$ with the following properties 
\begin{enumerate}
\item $\theta^2=\id$,
\item  $\theta(\De_+)= \De_-$ and, letting $c_\vep$ be the chamber corresponding to $\compB_\vep$ ($\vep=+,-$), we have $c_-=c_+^\theta$ and $c_- \opp c_+$,
\item $\theta$ preserves the twin-apartment $\Sigma(c_+,c_-)=(\Sigma_+(c_+,c_-),\Sigma_-(c_-,c_+))$ and 
permutes the types in each fiber of $p\colon \La\to \liediag$.
\end{enumerate}
\ele

\subsection{The Coxeter groups $W$ and $W^\theta$}
Here we recall some definitions and results from~\cite{BloHof2014a} and indicate how one can prove these in the current more general setting of $3$-spherical diagrams, which are not necessarily simply-laced.

Let $(W,\{s_j\colon j\in \hI\})$ be the Coxeter system of type $\La$ with $W=N/D$ and $s_j=\mu(u)D$ for $u\in U_{\alpha_j}$.
Let 
 $$\delta^\theta(W)=\{w\in W\mid \exists d_+\in \De_+\colon w=\delta_*(d_+,d_+^\theta)\}.$$
Let
 $$\Inv^\theta(W)=\{u\in W\mid u^\theta=u^{-1}\}.$$
and
 $$W(\theta)=\{w(w^{-1})^\theta \mid w\in W\}.$$

%
%

\ble\label{lem:tau does not preserve roots}
\label{lem:sws'<>w}
\label{lem:delta^theta(W)}
%
  $\Inv^\theta(W)=\delta^\theta(W).$
More precisely, given any $u\in \Inv^\theta(W)$ there exists a word $w\in W$ such that
 $w(w^{-1})^\theta$ is a reduced expression for $u$.
 [cf. Lemma 4.28 of ~\cite{BloHof2014a}]\label{lem:4.24BloHof}
\ele
This lemma is proved exactly as in~\cite{BloHof2014a} using Lemmas 4.24, 4.26, 4.27.
Note that the proof Lemma 4.24 relies on the fact that a connected spherical diagram does not admit an involution without fixed vertices or edges. As a consequence these results hold for any $3$-spherical diagram.

\bco\label{cor:every d in theta apartment}
Let $d\in \De_\vep$.
Then, there exists $(c,c^\theta)\in \De^\theta$ such that 
 $d\in \Sigma_+(c,c^\theta)$.
If the non-orientable Curtis-Tits amalgam $\amL$ is defined over $\FF_q$, there are at least $q^{\delta_+(c,d)}$ such chambers.
\eco
\pf
Let $u=\delta_*(d,d^\theta)$. We induct on $l(u)$.
If $u=1$, we are done.
Assume $l(u)>0$. From Lemma~\ref{lem:delta^theta(W)} we know that $u=w(w^{-1})^\theta$ is irreducible for some irreducible $w\in W$.
Pick any $s_i\in S$ such that $l(s_iw)=l(w)-1$ and let $e$ be a chamber $i$-adjacent to $d$.
Then, $\delta_*(e,e^\theta)=s_iw(w^{-1})^\theta s_{\theta(i)}$ is shorter by $2$, so by induction there exists
 $c$ such that $e\in \Sigma(c,c^\theta)_+$.
Moreover, calling $\pi$ the $i$-panel on $d$ we have
 $\proj^*_\pi(e^\theta)=d\in \Sigma(c,c^\theta)_+$. This in turn implies that if every panel has at least $q+1$ chambers, then there are at least $q$ choices for the pair $(e,e^\theta)$.
\qed
\bde
For each $u\in W^\delta$, define
\begin{align*}
\De_u^\theta&=\{d\in \De_+\colon \delta_*(d,d^\theta)=u\}.
\end{align*}
\ede

\bco\label{cor:comp transitive on De_u^theta}
The group $\compG$ acts transitively on the set $\De_u^\theta$ for each $u\in W^\delta$. Suppose that
 the orientable Curtis-Tits amalgam $\amL$ is defined over $\FF_q$, let  $d\in \De_u^\theta$ and suppose that $u=w(w^{-1})^\theta$ for some reduced $w\in W$.
Then, $|\Stab_\compG(d)|\ge q^{l(w)}$.
\eco
\bpf
It was proved in~\cite{BloHof2016} that $\compG$ is transitive on the set $\De_1^\theta$. 
The proof also applies in this case since all vertex groups are $\SL_2(q^e)$ acting on the panel as points of the projective line over $\FF_{q^e}$ for some $e\ge 1$.

Given $d\in \De_u^\theta$, let $c\in \De_1^\theta$ be such that $d\in \Sigma_+(c,c^\theta)$.
Then, $d$ is the unique chamber in $\Sigma_+(c,c^\theta)$ with $\delta_+(c,d)=w$.
The transitivity of $\compG$ on $\De_u^\theta$ now follows.

In addition, by Corollary~\ref{cor:every d in theta apartment}, given $d$, there are $q^{l(w)}$ such chambers $c$, hence at least $q^{l(w)}$ elements in $\compG$ fixing $d$.
\epf

\subsection{The group $G$ is a lattice in $\bhKM_+$}\label{subsec:lattice}
The action of the group $\compL$ on the positive building $\De_+$ turns it into a locally compact group with Haar measure $\mu$, and endowed  with a metric $f_+$ as in~\cite{CapRem2009a,CapRem2009}.  
Let $\bhKM_+$ be the completion of $\compL$ with respect to this metric.

\bpr\label{prop:G is a lattice of bhKM+}
The group $G$ is a lattice in the group $\bhKM_+$ for all $q\ge n$.
\epr
\bpf
We follow the proof idea in~\cite{BloHof2014a} and~\cite{GraMuh2008}.
The group $\bhKM_+$ is locally compact.
It follows from the fact that $U^+\cap \compG=\{1\}$ that $\compG$ acts discretely on $\De_+$ and {any two elements of $\compG$ have distance at least $2$ from each other in the metric $f_+$ defined there} (see also~\cite[\S 6]{CapRem2009a}).
We shall now show that  $\compG$ has finite covolume in $\bhKM_+$, that is $\compG\backslash \bhKM_+$ has finite volume.
To see this we use~\cite[Proposition 1.4.2]{Bou2000} and instead show that the sum
 \begin{align}\label{eqn:stabilizer series}
 \sum_{d\in \De_+}\frac{1}{|\Stab_\compG(d)|}
 \end{align}
 converges.
By Corollary~\ref{cor:comp transitive on De_u^theta} we see that if $d\in \De_+$ 
 then, $|\Stab_\compG(d)|\ge q^{l(w)}$.
 It follows that~\eqref{eqn:stabilizer series} is dominated by the Poincar\'e series of the Coxeter group $W$ evaluated at $t=q^{-1}$: 
 \begin{align}\label{eqn:evaluated poincare series}
p_{(W,S)}(q^{-1})=\sum_{l\in \NN}a_l(q^{-1})^l= \sum_{w\in W}\frac{1}{q^{l(w)}},
 \end{align}
 where $a_l$ is the number of elements in $W$ of length $l$ with respect to the generating set $S$.

Now since $S$ is a finite symmetric generating set for $W$, we have 
\begin{align*}
\rho_{(W,S)}^{-1}=\omega(W,S)=\limsup_{l\to\infty}\sqrt[l]{a_l},
\end{align*}
where $\omega(W,S)=\limsup_{l\to\infty}\sqrt[l]{a_l}$ is the growth rate of $(W,S)$ and $\rho_{(W,S)}$ is the radius of convergence of $p_{(W,S)}(t)$ as a power series over $\CC$.  
Thus if $q\ge \omega(W,S)$, then~\eqref{eqn:evaluated poincare series} converges.

We now wish to maximise $\omega(W,S)$ over all $3$-spherical $(W,S)$ with fixed $n=|S|$.  Recall the partial order $\preceq$ from ~\cite{Ter16}  on the set of Coxeter systems.
If $(W,S)$ and $(W',S')$ are two Coxeter systems with $S$ and $S'$ finite and Coxeter matrices $M$ and $M'$, then, we write 
 $(W,S)\preceq (W',S')$ whenever there is an injective map $\varphi\colon S\into S'$ such that for all $r,s\in S$ we have $m_{r,s}\le m'_{\varphi(r),\varphi(s)}$.
By~\cite[Theorem A]{Ter16} we then have $\omega(W,S)\le \omega(W',S')$.

Now suppose that  $(W',S')$ is a diagram dominating all $3$-spherical Coxeter systems with given $n=|S'|$ in the $\preceq$ order. Then~\eqref{eqn:evaluated poincare series} converges for all $q \ge \omega_{(W',S')}$.

For example, we can let $(W',S')$ be the Coxeter group whose diagram is the complete graph on $S'$ in which any two vertices are connected by a double edge.
To compute this series note that the set of spherical subsets of $\hI$ is  $\cF=\{J\sbe \hI\colon |W'_J|<\infty\}=\{J\sbe \hI\colon |J|\le 2\}$.
Then, by~\cite{St1968}, we have 
\begin{align}
\frac{1}{p_{(W',S')}(t^{-1})}=\sum_{J\in \cF}\frac{(-1)^{|J|}}{p_{(W'_J,S'_J)}(t)} = 1 -\frac{n}{p_1(t)}+{n\choose 2}\frac{1}{p_1(t)p_3(t)}, 
\end{align}
where, setting $p_m(t)=\sum_{i=0}^m t^i$, we have 
\begin{align*}
p_{(W'_J,S'_J)}(t) & = \begin{cases}
1 & \mbox{ if }J=\emptyset\\
p_1(t) & \mbox{ if }|J|=1\\
p_1(t)p_3(t) & \mbox{ if } |J|=2.
\end{cases}
\end{align*}
It follows that 
\begin{align*}
p_{(W',S')}(t)&=\frac{t^4}{t^4}\cdot \frac{2p_1(t^{-1})p_3(t^{-1})}{2p_1(t^{-1})p_3(t^{-1}) - 2np_3(t^{-1}) +n(n-1)}\\
&=\frac{2p_1(t)p_3(t)}{2p_1(t)p_3(t) - 2ntp_3(t) +t^4 n(n-1)}
\end{align*}
The denominator equals
\begin{align*}
2(1-(n-1)t)(1+t+t^2+t^3)+t^4n(n-1)
\end{align*}
When $|z|\le \frac{1}{n}$ one verifies that the last term has norm at most $\frac{1}{n^2}$, whereas the other terms together have norm at least $\frac{1}{n}$, noting that $n\ge 3$. Therefore the Poincar\'e series converges at $q^{-1}$ for $q\ge n$.
Numerical evidence suggests that one cannot do much better as, e.g.~for $n=1000$, the norm of the smallest complex root is about $0.001001...$.
 \epf

%




\section{Simplicity of Curtis-Tits groups}\label{sec:simplicity}
Combining Theorem~\ref{thm:OCT realization}, the Simplicity Theorem of~\cite{CapRem2009} and the observation that Curtis-Tits groups are perfect since they are generated by perfect subgroups, we have the following.
\bco\label{cor:OCT almost simple}
Suppose $\amL$ is an orientable Curtis-Tits amalgam over $\FF_q$ with connected diagram $\liediag$, which is $3$-spherical, but not spherical or affine. Then, its universal completion $\ucompL$ is almost simple.
\eco

By Theorem~\ref{thm:NOCT realization} (as well as~\cite{BloHofShp2017,BloHof2016}), any non-orientable Curtis-Tits amalgam $\amG$ has a completion inside the centralizer in a group of Kac-Moody type $\compL$ of a Cartan involution $\theta$. In contrast to what happens in the orientable case, we now have the following.

\bth\label{thm:non-orientable CT groups are not simple}
If $\compG$ is a non-orientable Curtis-Tits group over a finite field with irreducible non-spherical, non-affine diagram, then $\compG$ is acylindrically hyperbolic. In particular, it is not simple.
\eth
\bpf
By Remark 3.7 of~\cite{CapHum2015} it suffices to note the following.
The group $\compL$ is a group of Kac-Moody type with non-spherical, non-affine diagram over a finite field $\FF_q$.
The positive building associated to $\compL$,  denoted $\De_+$ is a proper CAT(0) space and $\Aut(\De_+)$ acts cocompactly on it.
Moreover, by~\cite[Theorem 1.1]{CapFuj2010} $\Aut(\De_+)$ contains rank 1 elements since it is of irreducible, non-spherical and non-affine type.
By~\cite[Corollary 3.6]{CapHum2015} any lattice of $\Aut(\De_+)$ is acylindrically hyperbolic and is therefore not simple.
Thus, the claim follows from Proposition~\ref{prop:G is a lattice of bhKM+}.
\epf
\subsection{Examples}\label{subsec:examples}
From the classification of Curtis-Tits amalgams it is clear that whenever the diagram is a tree, the amalgam is unique. In particular, the only Curtis-Tits groups with spherical diagram are the groups of Lie type.
The same holds for Curtis-Tits groups with affine diagram, other than $\tA_{n-1}$. 
In~\cite{BloHof2014a} all orientable and non-orientable Curtis-Tits groups with diagram $\tA_{n-1}$ were described in terms of matrix groups and it was shown that all of them have interesting quotients (see also~\cite{BloHofVdo2012}).

Here we will give some examples of quotients arising from non-orientable Curtis-Tits groups with non-spherical and non-affine diagram.
We are particularly interested in finite quotients.
In a subsequent study of such quotients we will use the finite presentations for Curtis-Tits groups arising from the Curtis-Tits amalgam and combine this with relations in the vein of~\cite{BarMar15}.
As a simple example consider the non-orientable Curtis-Tits amalgam over $\FF_2$ with the following diagram:
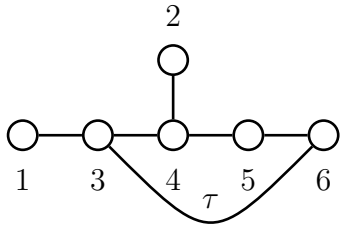
\begin{center}
 \begin{tikzpicture}[scale=.5]
\node [label=below:$\trin$] (tau) at (5,-1) {};
  \tikzstyle{every node} = [draw, line width = 1pt, shape=circle]
    \node  [label=below:$1$] (one) at (0,0) {};
   \node  [label=above:$2$] (two) at  (4,2)  {};
   \node  [label=below:$3$] (three) at (2,0) {};
   \node  [label=below:$4$] (four) at (4,0) {};
   \node  [label=below:$5$] (five) at (6,0) {};
   \node  [label=below:$6$] (six) at (8,0)  {};
\foreach \from/\to in {one/three, three/four, four/five, five/six, two/four}
\path[draw, line width = 1pt]  (\from) -- (\to);
 \path[draw, line width = 1pt] (three) .. controls (5,-3) .. (six);
  \end{tikzpicture}
\end{center}
Here the $\trin$ indicates that $\classmap$ sends the loop $\Lambda_1=\{3,4,5,6,3\}$ to $\trin$.
Let $\Chi=\{\amgrpX_i^+,\amgrpX_i^-\colon i=1,2,\ldots,6\}$ and, for each $i$, let $x_i^+$, $x_i^-$, and $n_i$ be defined by 
 $\amgrpX_i^+=\langle x_i^+\rangle$,  
$\amgrpX_i^-=\langle x_i^-\rangle$, and $n_i=x_i^+x_i^-x_i^+$.
Then, the universal completion has a presentation in terms of these generators. After adding the relation
\begin{align*}
(n_3n_4n_5n_6n_5n_4)^2=1,
\end{align*}
 GAP~\cite{GAP4.8.5} was able to identify the quotient as the exceptional group $E_6(2)$.
Since the orientable amalgam with the same diagram is simple by Corollary~\ref{cor:OCT almost simple}, this example confirms the existence of non-spherical, non-affine non-orientable Curtis-Tits groups.
Note that since $|\FF_2|<|S|=6$, this example is not covered by Proposition~\ref{prop:G is a lattice of bhKM+}.

We shall also study non-orientable Curtis-Tits amalgams with diagrams such as the following:
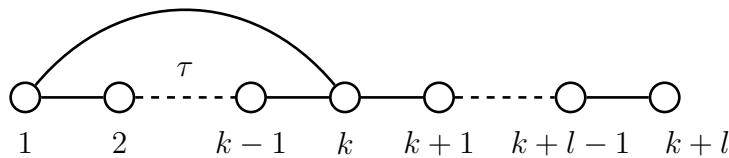
\begin{center}
 \begin{tikzpicture}[scale=.5]
\node [label=below:$\trin$] (tau) at (4.25,1.5) {};
  \tikzstyle{every node} = [draw, line width = 1pt, shape=circle]
   \node  [label= {[label distance=1ex] -90:\makebox(0,0){$1$}}]  (two) at (0,0) {};
   \node  [label= {[label distance=1ex] -90:\makebox(0,0){$2$}}]  (two) at (2.5,0) {};
   \node  [label= {[label distance=1ex] -90:\makebox(0,0){$k-1$}}]  (enminusone) at (6,0) {};
   \node  [label= {[label distance=1ex] -90:\makebox(0,0){$k$}}] (en) at (8.5,0) {};
   \node  [label= {[label distance=1ex] -90:\makebox(0,0){$k+1$}}]  (enplusone) at (11,0) {};
   \node  [label= {[label distance=1ex] -90:\makebox(0,0){$k+l-1$}}]  (enpluskayminusone) at (14.5,0) {};
   \node  [label= {[label distance=1ex] -90:\makebox(0,0)[l]{$k+l$}}]  (enpluskay) at (17,0) {};
\foreach \from/\to in {one/two, enminusone/en, en/enplusone, enpluskayminusone/enpluskay}
\path[draw, line width = 1pt]  (\from) -- (\to);
 \path[draw, line width = 1pt] (one) .. controls (2.5,3) and  (6,3)  .. (en);
 \foreach \from/\to in {two/enminusone, enplusone/enpluskayminusone}
\path[draw, dashed, line width = 1pt]  (\from) -- (\to);
 \end{tikzpicture}
\end{center}
Here the $\trin$ indicates that $\classmap$ sends the loop $\Lambda_1=\{1,2,\ldots,k,1\}$ to $\trin$.
Let $n_i$ ($i=1,2,\ldots,k+l$) be defined as above and consider the quotient of the universal completion of this amalgam over $\FF_q$ ($q$ even) over the normal closure of the element 
\begin{align*}
\left (n_1n_2\cdots n_{k-1}n_k n_{k-1}^{-1}\cdots n_2^{-1}\right)^2.
\end{align*}
Note here that the case $q=4$, $k=3$ and $l=1$, is covered by Proposition~\ref{prop:G is a lattice of bhKM+}, so that 
 the existence of infinitely many infinite-index normal subgroups is guaranteed via Theorem~\ref{thm:non-orientable CT groups are not simple}.
However, using GAP~\cite{GAP4.8.5} we were able to identify the following finite quotients.
For $k=3$ and $(q,l)= (2,1),(2,2),(2,3),(4,1),(4,2)$, the quotient is $\PSL_{k+l+1}(q)$.
For $k=4$,  $l=1$, and $q=2,4$, the quotient is $\POmega^+(2(k+l),q)$.


\begin{thebibliography}{1}
\bibitem{BarMar15}
M.~Barot and R.~J. Marsh.
\newblock Reflection group presentations arising from cluster algebras.
\newblock {\em Trans. Amer. Math. Soc.}, 367(3):1945--1967, 2015.



\bibitem{BloHofShp2017}
R.~J. Blok,  C.~G. Hoffman, and S.~V.~Shpectorov.
\newblock Classification of Curtis-Tits and Phan type amalgams with $3$-spherical diagram.
\newblock Submitted.

\bibitem{BloHof2016}
R.~J. Blok and C.~G. Hoffman.
\newblock Curtis-tits groups of simply-laced type.
\newblock To appear in J. Comb. Th. Ser. A., 2016.


\bibitem{BloHof2014a}
R.~J. Blok and C.~G. Hoffman.
\newblock Curtis--{T}its groups generalizing {K}ac--{M}oody groups of type
  {$\tilde{A}_{n-1}$}.
\newblock {\em J. Algebra}, 399:978--1012, 2014.

\bibitem{BloHof2014b}
R.~J. Blok and C.~G. Hoffman.
\newblock A classification of {C}urtis-{T}its amalgams.
\newblock In {\em Groups of exceptional type, {C}oxeter groups and related
  geometries}, volume~82 of {\em Springer Proc. Math. Stat.}, pages 1--26.
  Springer, New Delhi, 2014.


\bibitem{BloHofVdo2012}
R.~J. Blok, C.~G. Hoffman, and A.~Vdovina.
\newblock Expander graphs from {C}urtis-{T}its groups.
\newblock {\em J. Combin. Theory Ser. A}, 119(3):521--525, 2012.

\bibitem{Bou2000}
M.~Bourdon.
\newblock Sur les immeubles fuchsiens et leur type de quasi-isom\'etrie.
\newblock {\em Ergodic Theory Dynam. Systems}, 20(2):343--364, 2000.


\bibitem{Cap2007}
P.~E.~Caprace.
\newblock On 2-spherical {K}ac-{M}oody groups and their central extensions.
\newblock {\em Forum Math.}, 19(5):763--781, 2007.

\bibitem{CapFuj2010}
P.-E. Caprace and K.~Fujiwara.
\newblock Rank-one isometries of buildings and quasi-morphisms of {K}ac-{M}oody
  groups.
\newblock {\em Geom. Funct. Anal.}, 19(5):1296--1319, 2010.

\bibitem{CapHum2015}
P.-E. Caprace and D.~Hume.
\newblock Orthogonal forms of {K}ac-{M}oody groups are acylindrically
  hyperbolic.
\newblock {\em Ann. Inst. Fourier (Grenoble)}, 65(6):2613--2640, 2015.



%
\bibitem{CapRem2009a}
P.-E. Caprace and B.~R{\'e}my.
\newblock Groups with a root group datum.
\newblock {\em Innov. Incidence Geom.}, 9:5--77, 2009.


\bibitem{CapRem2009}
P.-E. Caprace and B.~R{\'e}my.
\newblock Simplicity and superrigidity of twin building lattices.
\newblock {\em Invent. Math.}, 176(1):169--221, 2009.


\bibitem{GAP4.8.5}
The GAP~Group.
\newblock {\em {GAP -- Groups, Algorithms, and Programming, Version 4.8.5}},
  2016.




\bibitem{GraMuh2008}
R.~Gramlich and B.~M\"uhlherr.
\newblock Lattices from involutions of kac-moody groups.
\newblock Report~5, Mathematisches Forschungsinstitut Oberwolfach, December
  2008.



\bibitem{Hoffman:2013aa}
C.~Hoffman and A.~Roberts.
\newblock On a quasi-{P}han theorem for orthogonal groups.
\newblock {\em Comm. Algebra}, 41(5):1589--1600, 2013.




\bibitem{Mu1999}
B.~M{\"u}hlherr.
\newblock Locally split and locally finite twin buildings of {$2$}-spherical
  type.
\newblock {\em J. Reine Angew. Math.}, 511:119--143, 1999.

\bibitem{MuRo1995}
B.~M{\"u}hlherr and M.~Ronan.
\newblock Local to global structure in twin buildings.
\newblock {\em Invent. Math.}, 122(1):71--81, 1995.



\bibitem{RemRon06}
B.~R{\'e}my and M.~Ronan.
\newblock Topological groups of {K}ac-{M}oody type, right-angled twinnings and
  their lattices.
\newblock {\em Comment. Math. Helv.}, 81(1):191--219, 2006.


\bibitem{St1968}
R.~Steinberg.
\newblock {\em Endomorphisms of linear algebraic groups}.
\newblock Memoirs of the American Mathematical Society, No. 80. American
  Mathematical Society, Providence, R.I., 1968.


\bibitem{Ter16}
T.~Terragni.
\newblock On the growth of a {C}oxeter group.
\newblock {\em Groups Geom. Dyn.}, 10(2):601--618, 2016.

\bibitem{Ti1992}
J.~Tits.
\newblock Twin buildings and groups of {K}ac-{M}oody type.
\newblock In {\em Groups, combinatorics \& geometry (Durham, 1990)}, volume 165
  of {\em London Math. Soc. Lecture Note Ser.}, pages 249--286. Cambridge Univ.
  Press, Cambridge, 1992.


\end{thebibliography}
\end{document}